\documentclass[12pt, letterpaper]{amsart}
\usepackage{amsmath}
\usepackage{amsthm}
\usepackage{amssymb}
\usepackage{amscd}
\usepackage{booktabs}
\usepackage{amsfonts}
\usepackage{graphicx}
\usepackage{graphics}
\usepackage{amsbsy}
\usepackage{epsfig}
\usepackage{amssymb,latexsym,mathrsfs}
\usepackage{array}
\usepackage[all]{xy}
\usepackage[colorlinks,linkcolor=blue,citecolor=red,linktocpage=true]{hyperref}
\usepackage{pgf,tikz}
\usepackage{mathrsfs}
\usetikzlibrary{arrows}
\usepackage{enumerate}
\usepackage[square,numbers]{natbib}

 \usepackage{setspace}

\hypersetup{citecolor=blue}
\DeclareGraphicsExtensions{.bmp,.png,.pdf,.jpg,.eps}
\graphicspath{{./figs/}}

\newtheorem {theorem} {Theorem}
\newtheorem {proposition} [theorem]{Proposition}

\newtheorem {remark} [theorem]{Remark}

 \usepackage{setspace}

\begin{document}

\title[On the equilateral pentagonal central configurations] {On the equilateral pentagonal central configurations}

\author[Alvarez-Ram\'{\i}rez]{M. Alvarez-Ram\'{\i}rez$^1$ }
 \address{$^1$ Departamento de Matem\'aticas, UAM--Iztapalapa\\
   09340 Iztapalapa,  Mexico City,  Mexico.  ORCID: 0000-0001-9187-1757}

\author[Gasul]{A. Gasull$^{2,3}$}

\author[Llibre]{J. Llibre$^{2}$}

\address{$^2$ Departament de Matem\`{a}tiques, Universitat Aut\`{o}noma de Barcelona, 08193 Bellaterra, Barcelona, Catalonia, Spain}
\address{$^3$ Centre de Recerca Matem\`{a}tica, Edifici Cc, Campus de Bellaterra, 08193 Cerdanyola del Vall\`{e}s (Barcelona),
Spain}

\email{mar@xanum.uam.mx, gasull@mat.uab.cat, jllibre@mat.uab.cat}

\subjclass[2010]{70F10 (70F15)}

\begin{abstract}
An equilateral pentagon is a polygon in the plane with five sides of equal length. In this paper we classify the central configurations of the $5$-body problem having the five bodies at the vertices of an equilateral pentagon with an axis of symmetry. We prove that there are two unique classes of such equilateral pentagons providing central configurations, one concave equilateral pentagon and one convex equilateral pentagon, the regular one. A key point of our proof is the use of rational parameterizations to transform the corresponding equations, which involve square roots, into polynomial equations.
\end{abstract}

\keywords{Central configuration, 5-body problem, equilateral pentagon}

\maketitle

\section{Introduction and statement of the result}

The Newtonian planar $5$-body problem describes the  dynamics of
five point particles  of positive masses  $m_i$ at positions
$\mathbf{q}_i\in \mathbb{R}^2$ moving according to the Newton's laws
under their mutual gravitational forces. The equations
of motion of this $5$-body problem are
\begin{equation*}\label{eq1}
m_i \ddot{\mathbf{q}}_i = -\sum_{j=1, j\neq i}^5 G m_im_j
\frac{\mathbf{q}_i - \mathbf{q}_j}{r_{ij}^3},    \hspace{1cm} 1\leq
i \leq  5,
\end{equation*}
where $r_{ij}= |\mathbf{q}_i-\mathbf{q}_j|$ is the mutual distances
between the masses $m_i$ and $m_j$, and $G$ is the gravitational
constant. We take conveniently the time unit so that $G=1$.

The configuration space is defined by
$$
{\mathcal E}  = \{{\mathbf q}=({\mathbf q}_1,\dots, {\mathbf q}_5)\in (\mathbb{R}^2)^5: {\mathbf q}_i\neq {\mathbf q}_j,   \quad  i\neq j\}.
$$
The configuration $\mathbf{q}=(\mathbf{q}_1,\dots, \mathbf{q}_5)$ is
called {\em central} if the position vector of each body with
respect to the center of mass is proportional to the corresponding
acceleration vector. In other words, if there exists a positive
constant  $\lambda$ such that
\begin{equation*}\label{eq_cc1}
\ddot{\mathbf{q}}_i  = \lambda ({\mathbf q}_i-{\mathbf c_m}), \qquad i=1,\dots, 5,
\end{equation*}
where ${\mathbf c}_m=(m_1{\mathbf q}_1+\dots + m_5 {\mathbf q}_5)/M$
and $M=m_1+\dots + m_5$,  being ${\mathbf c}_m$ and $M$ the center
of mass of the five bodies and the total mass, respectively. Hence a
given configuration $(\mathbf{q}_1,\dots, \mathbf{q}_5)\in {\mathcal
E}$ of the $5$-body problem with positive masses $m_1,\dots, m_5$,
is central if there exists a $\lambda$ such that
$(\lambda,\mathbf{q}_1,\dots, \mathbf{q}_5)$ is a solution of the
system
\begin{equation}\label{cc_2}
\sum_{j=1, j\neq i}^5m_j \frac{\mathbf{q}_i -
\mathbf{q}_j}{r_{ij}^3}=\lambda (\mathbf{q}_i - \mathbf{c}_m) ,
\hspace{1cm} 1\leq i \leq  5.
\end{equation}

A central configuration is {\em convex} if no body belongs to  the
convex hull of the other four bodies; otherwise it is called {\em
concave}. A planar central configuration is called a {\em relative
equilibrium} when they become equilibrium solutions in a rotating
coordinate system \cite{Meyer}.

Since equations \eqref{cc_2} are invariant under rotations,
translations and dilations,  when we consider the number of central
configurations, this will be restricted to count the classes of
central configurations modulo these mentioned transformations.

The central configurations are of special importance in  Celestial
Mechanics for several reasons. For instance, central configurations
are the initial conditions for the homographic orbits of the
$n$-body problem. Central configurations play an important role in
the description of the topology of the integral manifolds in the
$n$-body problem. Moreover, in the planar case the central
configurations are initial positions for periodic solutions. For
more information on this subject, recent advances and open
questions, the reader is addressed to \cite{moeckel2015} and
references therein.

In this paper we will investigate some central configurations  of
the $5$-body problem, for which there are few known results. The
first results concerning this issue take us back decades ago to the
work due to Williams \cite{Williams1938}, who settled  necessary and
sufficient conditions for  any plane central configuration of five
bodies. In what follows we offer a non-exhaustive list of some
interesting works concerning this topic, which have been published
more recently.

Albouy and Kaloshin \cite{AlbouyKaloshin} proved that for a choice
of  five positive masses in the complement of a codimension-two
algebraic variety in the mass space, there are only a finite number
of equivalence classes of central configurations of the Newtonian
5-body problem. In \cite{Chen2018} Chen and  Hsiao provided
necessary conditions for strictly convex central configurations of
the planar 5-body problem.

In the last times the interest in stacked central configurations
has grown a lot, that is,  central configurations in which some
subset of three or more masses also forms a central configuration.
This concept was introduced by  Hampton  \cite{Hampton2005}, who was
arguably the first to find stacked central configurations in the
5-body problem, where  two bodies can be removed and the remaining
three bodies are already in a central configuration. After, several
papers have been published showing the existence of other stacked
central configurations in the planar 5-body problem;  see, among
others, \cite{Lino2017}, \cite{Lino2019, Mello2013,  GideaLlibre2010, LMC,
Mello2008}.

Other studies have focused on restricting the problem to a
particular shape, in \cite{LlibreValls2015} it was proved that the
unique co-circular central configuration in the planar 5-body
problem is the regular 5-gon with equal masses, while in
\cite{Mello2008} was proved the existence of three  families of
planar central configurations where three bodies are at the vertices
of an equilateral triangle and the other two bodies are on a
perpendicular bisector. Later on in \cite{Shoaib2017} was studied
the central configuration in a symmetric $5$-body problem with three
masses on an axis of symmetry and two other masses outside this
axis, placed in symmetric positions. A complete  classification  of
the isolated central configurations of the planar $5$-body problem
with equal masses was given in \cite{Santoprete2009}. Recently in
\cite{Lino2021} were studied  the central configurations of the
planar $5$-body problem having four bodies at the vertices of a
rhombus.

\begin{figure}[h]
\begin{center}
\includegraphics[width=0.7\textwidth]{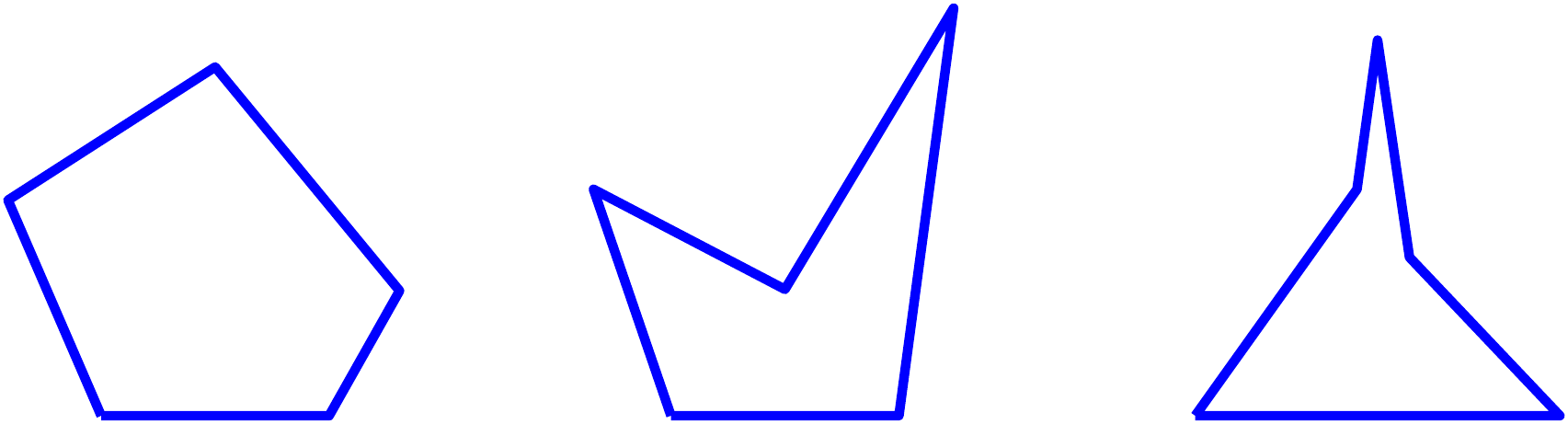}
\end{center}
\caption{Pentagonal configurations: a convex one in the left hand picture  and two different type of concave ones} \label{graf1}
\end{figure}

For the reader's convenience, we summarize here some basic facts
about the pentagon, which  is a polygon  with five sides and five
angles. It is said that it is {\em convex} if all its vertices are
pointing outwards, otherwise  it is {\em concave}; see Figure
\ref{graf1}. A pentagon with five sides of equal length is named
{\em equilateral}. Moreover, a pentagon is called {\em regular} when
all the sides are equal in length, and five angles are of equal
measures. If the pentagon does not have equal side length and angle
measure, then it is called {\em irregular}. The regular pentagon is
unique, up to similarity transformations, because it is equilateral
and its five angles are equal.

The goal in this paper is to characterize the equilateral pentagonal
central configurations with an axis of symmetry for the planar
$5$-body problem whose five positive masses are at the vertices of
an equilateral pentagon.

\begin{figure}[h]
\begin{center}
    \includegraphics[width=0.7\textwidth]{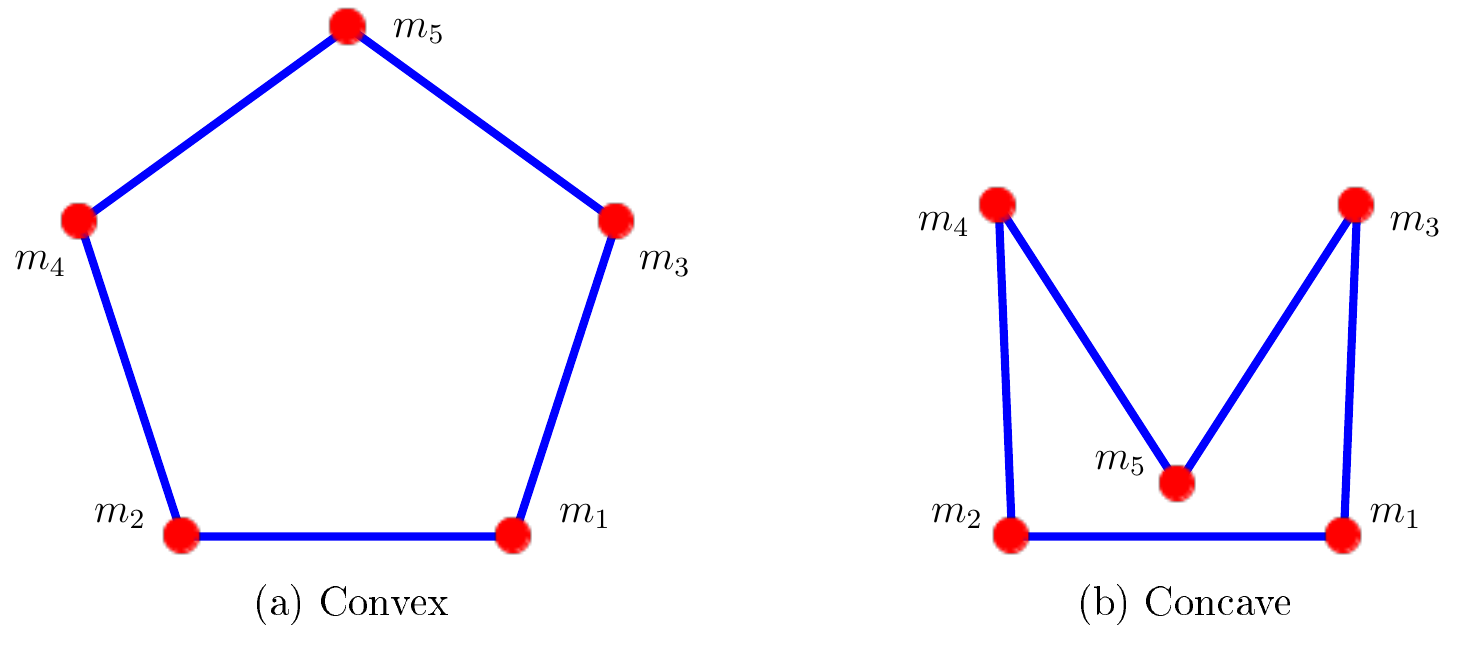}
\end{center}
\caption{The two equilateral pentagonal central configurations with an axis of symmetry of the $5$-body problem.}
\end{figure}

Then our main result is the following one, which will be proved in the next section.

\begin{theorem}\label{t1}
There are two classes of equilateral pentagonal central configuration having an axis of symmetry for the $5$-body problem.
\begin{itemize}
\item[(a)] The convex regular pentagon with equal five masses, see Figure $2$(a).

\smallskip

\item[(b)] The equilateral concave  pentagon with the masses
normalized, i.e. $\sum_{i=1}^5 m_i=1$, equal to $m_1=m_2\approx 0.0922539749$, $m_3=m_4\approx 0.3860948766$ and $m_5\approx 0.04330242730.$ In  Figure $2$(b) we show a representative of its class where the bodies of masses $m_1$ and $m_2$ are fixed at $(0,1/2)$ and $(0,-1/2),$ respectively. The other bodies of masses $m_3,m_4$ and $m_5$ are located at $(x_3,y_3),$ $(-x_3,y_3)$ and $(0,y_5),$ respectively, where $x_3\approx0.5402091568,$ $y_3\approx0.9991912848$ and $y_5\approx0.1576604970.$ See Remark $\ref{re:1}$ for some comments of how these values are obtained.
\end{itemize}
\end{theorem}

\begin{remark}\label{re:1}
From the proof of Theorem $\ref{t1}$ and the results of Appendix \ref{ap2} we obtain a value $t=t^*\approx 0.7332148086$ which is the smallest root of the quadratic polynomial $t^2-u^*t+1=0,$ where $u=u^*\approx 2.0970716051$ is the unique root in the interval $[205/100,210/100]$ of the polynomial of degree $60$
with integer coefficients, $R_{60}(u),$ given in Appendix \ref{ap2}. This polynomial is constructed from a reciprocal polynomial of degree $120$ and also  with integer coefficients that appears in \eqref{eq:pt}. Then $y_5=(1-(t^*)^2)/(4t^*)$ and all the other values in the theorem, $x_3,y_3,m_j, j=1,2,\ldots,5,$ can be obtained from this $t^*$ by elementary computations: sums, subtractions, multiplications, divisions and square roots.  Recall that, although we do not know the exact value of $u^*,$ the classical use of Sturm sequences allows to obtain explicit intervals, with rational endpoints, containing $u^*$ and with arbitrarily small length.
\end{remark}

At this stage the reader should be warned that there is a previous
work by  Perko and Walter \cite{Perko}, who showed that $n$ equal
masses at the vertices of a regular polygon, for $n\geq 4$, forms a
central configuration if and only if the masses are equal. Therefore
it was known that the regular pentagon with equal masses is a
central configuration for the $5$-body problem, but it was unknown
that it is the unique equilateral convex pentagonal central
configuration with an axis of symmetry.

One of the key points of our approach is the use of rational
parameterizations to eliminate some of the square roots that appear
in the equations governing the central configurations, converting in
this way these equations into polynomial ones. Then these equations
can be treated analytically by using some classical tools, like for
instance the Sturm sequences or the computation of resultants.

\section{Preliminaries}

Central configurations are invariant under composition of translations,
 rotations, and scaling through its center of mass, hence without loss of
 generality we can assume that the position of the masses $m_i> 0$ for $i=1\ldots,5$
 at the vertices of an equilateral pentagon with an axis of symmetry are $p_k=(x_k,y_k)$
 for $k=1, 2, 3, 4, 5$, where $y_1=0$, $x_2=-x_1$, $y_2=0$, $x_4=-x_3$, $y_4=y_3$ and $x_5=0$.
 Note that we can assume that $x_1>0$, $x_3>0$, $y_3> 0$, $y_5> 0$ and $y_3\ne y_5$, because we
 want that the points $(x_i,y_i)$ be the vertices of a pentagon.

Next we will obtain the coordinates for the equilateral pentagon
vertices. By a suitable scaling we may assume that  $r_{12}=1$, so  $x_1= 1/2$.
Now we substitute this value into the equation $r_{13}=1$,
obtaining that $y_3^2= (3 + 4 x_3-4 x_3^2)/4$. These values
replaced in the equation $r_{35}=1$ provides $x_3$ in terms of
$y_5$,
\begin{equation}\label{eq:psi}
x_3=\Psi^\pm(y_5)=\frac14 \pm\frac{y_5}2\Phi(y_5),\quad\mbox{where}\quad \Phi(y)=\sqrt{\frac{15-4y^2}{1+4y^2}}.
\end{equation}

Since we are interested in central configurations, with an axis of symmetry, modulus rotations and homothetic transformations,
without loss of generality, we have the next result.
\begin{proposition}\label{p1}
Suppose that $\mathbf{q}_1, \dots, \mathbf{q}_5$ form an equilateral pentagon,  with  $\mathbf{q}_1$ and $\mathbf{q}_2$   fixed  on the $x$-axis and  with the $y$-axis as an axis of symmetry. Then
\[
(x_1,y_1)=(0,1/2),\quad (x_2,y_2)=(-1/2,0), \quad (x_3,y_3), \quad (-x_3,y_3)\quad \mbox{and}\quad (0,y_5),
\]
where
\begin{enumerate}
\item[{(a)}] either
$
x_3=\Psi^+(y_5)=\dfrac14 +\dfrac{y_5}2\Phi(y_5),$  \quad $y_3=\dfrac{y_5}2 +\dfrac{1}4\Phi(y_5), \quad\mbox{and}\quad y_5\in\big(0,\sqrt{15}/2\big),
$
\item[{(b)}] or
$
x_3=\Psi^-(y_5)=\dfrac14 -\dfrac{y_5}2\Phi(y_5) ,$ \quad $y_3=\dfrac{y_5}2 -\dfrac{1}4\Phi(y_5),  \quad\mbox{and}\quad y_5\in\big(1+\sqrt3/2,\sqrt{15}/2\big).
$
\end{enumerate}
\end{proposition}

\begin{figure}[h]
    \begin{center}
        \includegraphics[width=0.83\textwidth]{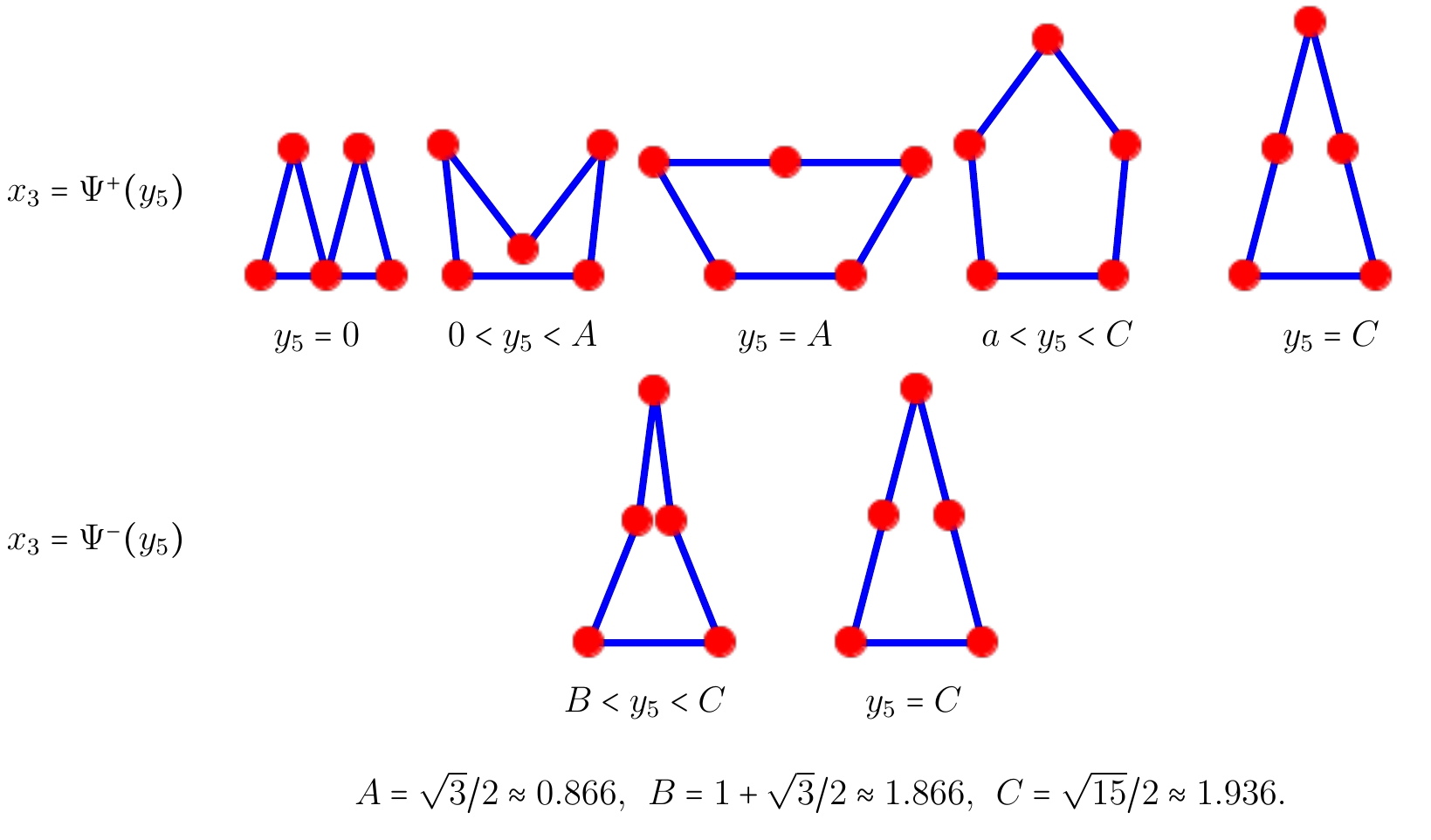}
    \end{center}
    \caption{Gallery of possible pentagon equilateral configurations, according whether $x_3=\Psi^+(y_5)$ or $x_3=\Psi^-(y_5)$ and the value of $y_5$, together with the boundary cases, that are no more pentagons.}
    \label{figy5}
\end{figure}

The geometrically distinct equilateral pentagon are shown in Figure \ref{figy5}. It follows from
the case (a) of Proposition \ref{p1} that  the equilateral pentagon is concave when $y_5\in (0,\sqrt{3}/2)$,
and convex if $y_5\in (\sqrt{3}/2,\sqrt{15}/2)$.
While in case (b) of Proposition \ref{p1} the equilateral pentagon is always concave.

In the next section we will see that the only values of $y_5$ that give rise to central configurations will be $x_3=\Psi^+(y_5)$, where $y_5=(\sqrt{5+2\sqrt5})/2\approx 1.539$ is associated to the regular pentagon, while $y_5\approx0.1576605$  gives a concave central configuration.

\section{Proof of Theorem $\ref{t1}$}

By a suitable scaling we may assume that  $m_1 + m_2 + m_3 + m_4 + m_5 =1$. Then the center of mass of the five bodies is
$$
{\mathbf c}_m = (x_m,y_m)=\left(\frac{m_1-m_2+2(m_3-m_4)x_3}{2},  \frac{(m_3+m_4)\sqrt{3+4x_3-4x_3^2}+2m_5y_5}{2}\right).
$$

Since we are studying equilateral pentagons we can assume that
$r_{12}=r_{13}=r_{35}=r_{45}=r_{24}=1$. Then the others mutual
distances are
\begin{equation}\label{con_dist}
r_{14}=r_{2,3}=\sqrt{1+2x_3},\qquad
r_{1,5}=r_{2,5}=\sqrt{y_5^2+1/4},\qquad r_{3,4}=2x_3.
\end{equation}
From \eqref{cc_2} we obtain the ten equations for the central configurations of the 5-body problem in the plane:
\begin{equation}\label{eqs}
\begin{array}{l}
e_j= \displaystyle{\sum_{j=1, j\neq i}^5\frac{m_j (x_i-x_j)}{r_{i,j}^3} - \lambda (x_j- x_m)}=0, \quad 1 \leq j\leq 5, \vspace{0.3cm}\\
e_{j+5} = \displaystyle{\sum_{j=1, j\neq i}^5\frac{m_j
(y_i-y_j)}{r_{i,j}^3} - \lambda (y_j- y_m)}=0, \quad 1 \leq j\leq 5.
\end{array}
\end{equation}

Substituting into \eqref{eqs} the values and expressions  of the mutual distances and taking $m_5=1-m_1-m_2-m_3-m_4$, it is seen that
\begin{equation}\label{xya}
e_8 - e_9 =-\frac{(m_1-m_2)\sqrt{3+4x_3-4x_3^2} (-1+\sqrt{1+2x_3}+2x_3\sqrt{1+2x_3}\:)}{2(1+2x_3)^{3/2}}=0.
\end{equation}
A straightforward computation shows that $\sqrt{3+4x_3-4x_3^2}
(-1+\sqrt{1+2x_3}+2x_3\sqrt{1+2x_3})=0$ for $x_3=-1/2,0,3/2$.
However, any of these values is good, because $x_3>0$ and $x_3=3/2$
implies that $r_{1,3}=\sqrt{1+y_3^2}$, but this is impossible
because we have assumed that  $y_3\ne0$ and $r_{1,3}=1$. So
$m_2=m_1$.

Since $e_3 + e_4 =(m_3-m_4)(1+8\lambda x_3^3)/(4x_3^2)=0$. It follows that either $m_4=m_3$, or $\lambda = -1/(8x_3^3)$. In this last case we obtain
$$
e_6 - e_7 = \frac{(m_3-m_4)\sqrt{3+4x_3-4x_3^2} (-1+\sqrt{1+2x_3}+2x_3\sqrt{1+2x_3}\:)}{2(1+2x_3)^{3/2}}=0.
$$
Hence as in \eqref{xya} we have that $m_4=m_3$. Therefore we do
not need to consider $\lambda = -1/(8x_3^3)$,  and in what follows
we consider that $m_4=m_3$, such that,   $e_6-e_7=0$.

In summary, we have that $e_1+e_2=0$, $e_5=0$ and $e_6-e_7=0$. Hence
we conclude that from the ten  equations \eqref{eqs} only  $e_1,
e_3, e_6, e_8, e_{10}$ remain independent. These equations are
$$
\begin{array}{rl}
f_1 =& -\dfrac{\lambda}{2}-\dfrac{4}{E_3}  + m_1 \left(-1 +\dfrac{8}{E_3}\right) + m_3 \left(-\dfrac{1}{2}+x_3-\dfrac{1}{2E_1}
+\dfrac{8}{E_3}\right),\vspace{0.2cm}\\
f_2 =& -(1+\lambda) x_3+\dfrac{1}{4}m_3 \left(8 -\dfrac{1}{x_3^3} \right)x_3+m_1 \left(\dfrac{1}{2}+x_3-\dfrac{1}{2E_1} \right),\vspace{0.2cm}\\
f_3 =& y_5  \left( \lambda +\dfrac{8}{E_3} \right)  + m_1 \left( -2\lambda y_5 -\dfrac{16y_5}{E_3}\right)\vspace{0.2cm}\\
&  +m_3 \left( \dfrac{1}{2} E_2\left( 1+\dfrac{1}{(1+2x_3)E_1}\right)+ \lambda \Big(E_2 -2y_5\Big)
-\dfrac{16y_5}{E_3} \right),\vspace{0.2cm}\\
f_4 =&  -\dfrac{1}{2} (1+\lambda)  \Big(E_2 -2y_5\Big) + (1+\lambda )m_3 \Big(E_2 -2y_5\Big)\vspace{0.2cm}\\
&+\; m_1 \left(\dfrac{1}{2}E_2\left( 1- \dfrac{1}{(1+2x_3)E_1}\right) -2 (1+\lambda) y_5  \right),\vspace{0.2cm}\\
f_5 =& (1+\lambda) m_3 \left( E_2 -2y_5\right) + m_1 \left( -2\lambda y_5 -\dfrac{16y_5}{E_3}\right),
\end{array}
$$
where $E_1= \sqrt{1 + 2 x_3}$, $E_2= \sqrt{3 + 4 x_3 - 4 x_3^2}$ and  $E_3= (1 + 4 y_5^2)^{3/2}$.

Solving $f_2=0$ and $f_5=0$, we obtain the following expressions for
$m_1$ and $m_3$
\begin{equation}\label{eq:m1m3}
\begin{array}{c}
m_1 = \big(2 E_1 E_3 (1 + \lambda)^2 x_3^3 (E_2 - 2 y_5))/m,\vspace{0.3cm}\\
 m_3=
\big( 4 E_1 (1 +\lambda) (8 + E_3 \lambda) x_3^3 y_5)\big)/m,
\end{array}
\end{equation}
where
\begin{multline*}
m=2 {   E_1}  \left( {   E_2} {   E_3} \lambda+2 \lambda {   y_5}
 {   E_3}+{   E_2} {   E_3}-2 {   E_3} {   y_5}+32 {   y_5}
\right) {{   x_3}}^{3}\\+{   E_3}  \left( 1+\lambda \right)  \left( {
       E_2}-2 {   y_5} \right)  \left( {   E_1}-1 \right) {{   x_3}}^{2}-
{   E_1} {   y_5}  \left( \lambda {   E_3}+8 \right).
\end{multline*}

We substitute the values of $m_1$ and $m_3$ into the  equations
$f_1=0$, $f_3=0$ and $f_4=0$, and taking  only the numerators of these three equations because
the denominators do not vanish, the former system, reduce to
\begin{align*}
g_1=&8 {  E_1} {  E_3}  \left( 1+\lambda \right)  \left( \lambda {
      E_3}+8 \right) {  y_5} {{  x_3}}^{4}+ \Big(  \big( 16 {  E_1
} {{  E_3}}^{2}\lambda-4 {\lambda}^{2}{{  E_3}}^{2}+8 {  E_1} {
    {  E_3}}^{2}-192 \lambda {  E_1} {  E_3}-4 \lambda {{  E_3}}^
{2}\\&-64 {  E_1} {  E_3}+512 {  E_1} \lambda-32 \lambda {  E_3
}-32 {  E_3} \big) {  y_5}-2 {  E_1} {  E_2} {  E_3}
\left( 1+\lambda \right)  \left( 3 \lambda {  E_3}+2 {  E_3}-16
 \lambda-8 \right)  \Big) {{  x_3}}^{3}\\&+ \left( 2 {  E_3}
\left( 1+\lambda \right)  \left( \lambda {  E_3}+8 \right)  \left( {  E_1}-1 \right) {  y_5}-{  E_2} {  E_3}  \left( 1+\lambda
\right)  \left( \lambda {  E_3}+8 \right)  \left( {  E_1}-1
\right)  \right) {{  x_3}}^{2}\\&+{  E_1}  \left( \lambda {  E_3}+8
\right) ^{2}{  y_5},\\
g_3=& (8 + E_3 \lambda) y_5 (L_2\lambda-L_1),\\
g_4=& (1
+ \lambda) (2 y_5-E_2) (L_2\lambda-L_1),
\end{align*}
where
\begin{align*}
L_1=&-2 \left( 1+2{ x_3} \right)  \left( 2{ E_1}{ E_3}{{
         x_3}}^{3}-{ E_1}{ E_3}{{ x_3}}^{2}+{{ x_3}}^{2}{ E_3}
-4{ E_1} \right) { y_5}-{ E_2}{ E_3}{{ x_3}}^{2} \left( 2{ x_3}{
E_1}+{ E_1}-1 \right),
\vspace{0.2cm}  \\
L_2=&  \left( 1+2{ x_3} \right)  \left( 12{ E_1}{ E_3}{{ x_3
}}^{3}-2{ E_1}{ E_3}{{ x_3}}^{2}-64{{ x_3}}^{3}{ E_1 }+2{{
x_3}}^{2}{ E_3}-{ E_1}{ E_3} \right) { y_5}\\&+{
    E_2}{ E_3}{{ x_3}}^{2} \left( 2{ x_3}{ E_1}+{ E_1}-1
\right),
\end{align*}

Computing $\lambda$ from equation $g_3=0$ we obtain the two solutions
\begin{equation}\label{eq:lambda}
\lambda_1=  -\frac{8}{E_3},  \qquad \lambda_2= \frac{L_1}{L_2}.
\end{equation}
 The solution $\lambda=\lambda_1$ is not suitable because then $m_3=0$.
Substituting the  solution $\lambda=\lambda_2$ in the equations $g_1=0$ and
$g_4=0$, we get that $g_4\equiv 0$, and  the equation $g_1=0$
reduces to
$$
\bar h_1= (2x_3^3 y_5 E_1E_3  (E_3-8)^2 (2 x_3-1) (1 + 2 x_3 + 4
x_3^2)
  h_1=0,
$$
where
\begin{align*}
  h_1=&-4 \left( 1+2{ x_3} \right) ^{2} \big( 4{{ E_1}}^{2}
\left( { E_3}-16 \right) {{ x_3}}^{4}-8{ E_1 } \left( { E_1}{
E_3}-4{ E_1}-4 \right) {{ x_3}}^{3}-{
     E_3} \left( { E_1}-1 \right) ^{2}{{ x_3}}^{2}\\&+{{ E_1}}^{2}
\left( { E_3}-8 \right)
\big) {{ y_5}}^{2}-2{ E_2} \left( 1+2{
     x_3} \right)  \big(
16{{ E_1}}^{2} \left( { E_3}-4 \right) {{ x_3}}^{4}+4{ E_1
} \left( 3{ E_1}{ E_3}-8{ E_1}-{ E_3}-8 \right) {{
        x_3}}^{3}\\&+2{ E_3} \left( { E_1}-1 \right) ^{2}{{ x_3}}^{2}-2
{{ E_1}}^{2} \left( { E_3}-8 \right) { x_3}-{{ E_1}}^{2}
\left( { E_3}-8 \right)
 \big) { y_5}+{{ E_2}}^{2}{
    E_3}{{ x_3}}^{2} \left( 2{ x_3}{ E_1}+{ E_1}-1 \right) ^{
    2}.
\end{align*}
Notice that $2x_3^3 y_5 E_1 E_3  (1 + 2 x_3 + 4 x_3^2)$ does not
vanish because $x_3\in (0,1)$,  $y_5>0$ and $E_1=r_{1,4}>0.$

At this step we shall prove that  condition $(E_3-8)(2 x_3-1)=0$ implies that $m_5=0.$ This
is so, because
\begin{align*}
m_5=&1-(m_1+m_2+m_3+m_4)=1-2(m_1+m_3)\\=&\frac{2 { E_1} { y_5} {{
x_3}}^{3}{ E_2} { E_3}  \left( 2 {  x_3}-1 \right) \left( 4 {{
x_3}}^{2}+2 { x_3}+1 \right)
 \left( { E_3}-8 \right)
}{L_2^2}\\&\times\Big( - \left( 1+2 { x_3} \right)  \Big[ 8 {{
E_1}}^{2} \left( { E_3}-16 \right) {{ x_3}}^{4}+8 { E_1 }  \left( {
E_3}-8 \right)  \left( { E_1}+1 \right) {{ x_3}}^{3 }+2 { E_3}
\left( { E_1}-1 \right) ^{2}{{ x_3}}^{2}\\&-2 {{ E_1}}^{2} \left( {
E_3}-16 \right) { x_3}-{ E_1}  \left( { E_1 } { E_3}-16 { E_1}+{
E_3} \right) \Big] { y_5}+{ E_2} { E_3} {{ x_3}}^{2} \left( 2 { E_1}
{ x_3}+{ E_1}-1 \right) ^{2} \Big),
\end{align*}
where we have used the expressions \eqref{eq:m1m3} for $m_1,m_3$
and substituted $\lambda=\lambda_2$ where $\lambda_2$ is given in
equation \eqref{eq:lambda}. In short we have proved that the equations
$\bar h_1=0$ and  $ h_1=0$ are equivalent.

From Proposition \ref{p1} we have that
$$
x_3=\frac14 \pm\frac{y_5}2\sqrt{\frac{15-4y_5^2}{1+4y_5^2}},
$$
or equivalently,
\begin{equation}\label{eq:h2}
h_2= (1 - 4 x_3)^2 (1 + 4 y_5^2) + 4 y_5^2 (4y_5^2-15)=0.
\end{equation}

Hence the central configurations are the solutions of the simultaneous solution of both equations $h_1=0$ and $h_2=0,$ with the two unknowns $x_3$ and $y_5.$  Indeed it provides positive masses $m_j.$ In order to avoid
the square roots which appear in  $E_1,E_2$ and $E_3$ in $h_1$, we do
a change of variables such that the expressions  appearing inside
each square root are equal to some new squared expressions. These
changes of variables are given by  the so called rational
parameterizations and correspond to parameterizations of planar
algebraic curves given by rational functions.  Due to the famous
Cayley-Riemann's Theorem \cite{Abh,Abh2} they exist if and only if
the corresponding surfaces have  genus zero. There
are effective methods to find one of these parameterizations see for instance \cite[Chap. 4
\& 5]{Sen}. In fact, many programs of symbolic calculus have
implemented some  methods and algorithms for obtaining them.
For more information on this subject the reader is addressed  to  \cite{GasLazTor} and references therein, where there are
several examples of applications of this approach.

In our case, for instance, we have that $E_1= \sqrt{1 + 2 x_3}$,
$E_2= \sqrt{3 + 4 x_3 - 4 x_3^2}.$  Hence if we write
$x_3=(u^2-1)/2$ we get that $E_1=\sqrt{u^2}.$ Then
\[
3 + 4 x_3 - 4 x_3^2\Big|_{x_3=(u^2-1)/2}=u^2(4-u^2).
\]
Consider now the algebraic curve $F(u,v)=u^2(4-u^2)-v^2=0.$ It has
genus 0, and by the Cayley-Riemann's theorem it admits a rational
parameterization. For example, for all $s$,
\[
F\Big({\frac {4\left( 2s-1 \right) s}{5{s}^{2}-4s+1}},{\frac {8s
\left( s-1 \right)  \left( 3s-1 \right)  \left( 2s-1
 \right) }{ \left( 5{s}^{2}-4s+1 \right) ^{2}}}
\Big)=0.
\]
As a consequence, by taking
\[
x_3=u^2(4-u^2)\Big|_{u={\frac {4\left( 2s-1 \right)
s}{5{s}^{2}-4s+1}}}=\frac{\left(3 s^2-1\right) \left(13 s^2-8
s+1\right)}{2 \left(5 s^2-4 s+1\right)^2}
\]
we obtain that
\[
E_2=\sqrt{v^2}=\sqrt{{\frac {\big(8s \left( s-1 \right)  \left( 3s-1
\right)  \left( 2s-1
 \right) \big)^2}{ \left( 5{s}^{2}-4s+1 \right) ^{4}}}}.
\]
The rational parameterization of $E_3=(1+4y_5^2)^{3/2}$ is much
simpler and can be obtained similarly. It suffices to consider the
algebraic curve of genus 0, $G(y_5,w)=1+4y_5^2-w^2=0.$  A good
parameterization for $y_5$ is  $y_5=(1-t^2)/(4t).$

We are interested on values $y_5\in (0,\sqrt{15}/2)$ and
$x_3\in(0,1)$. In short, we do the change of variables $(x_3,y_5)\to
(s,t)$ where
$$
x_3= \frac{\left(3 s^2-1\right) \left(13 s^2-8 s+1\right)}{2 \left(5
s^2-4 s+1\right)^2},\quad y_5= \frac{1-t^2}{4 t}.
$$
Then we have that $t$ varies in the interval $T:=(4 -\sqrt{15},1)$
and similarly, the values of $x_3\in(0,1)$ are covered for instance
for $s\in S:=(\sqrt3/3, (6+\sqrt3)/11).$ Hence,
$$
E_1=\frac{4 s (2 s-1)}{5 s^2-4 s+1},\quad E_2=\frac{-8 s
(s-1)(3s-1)(2s-1)}{\left(5 s^2-4 s+1\right)^2}, \quad
E_3=\frac{\left(t^2+1\right)^3}{8 t^3}.
$$

In the variables $(s,t)$ the two equations $ h_1=0$ and $h_2=0$
become
\begin{equation}\label{r1}
\begin{array}{rl}
h_1=&  \dfrac{8s^2 (2s-1)^2}{\left(5 s^2-4 s+1\right)^2 t^3} H_1=0,  \vspace{0.2cm} \\
h_2=&  \dfrac{1}{16 \left(5 s^2-4 s+1\right)^4 t^4} H_2=0,
\end{array}
\end{equation}
respectively, where
\begin{align*}
H_2=&\left( 5 {s}^{2}-4 s+1 \right) ^{4}\big({t}^{8}+1\big)-4 \left(
47 {s}^{4}- 152 {s}^{3}+150 {s}^{2}-56 s+7 \right)\\&\times
\left( 153 {s}^{4}-168 {s }^{3}+58 {s}^{2}-8 s+1 \right)
\big({t}^{6}+t^2\big) +\big( 101222 {s}^{8}- 258784 {s}^{7}+326904
{s}^{6}\\&-286240 {s}^{5}+183428 {s}^{4}-79776  {s}^{3}+21432
{s}^{2}-3168 s+198 \big) {t}^{4}
\end{align*}
and $H_1$ ia a huge polynomial of total degree 34 and with
$\operatorname{deg}_s(H_1)=24$ and $\operatorname{deg}_t(H_1)=10,$
which is given in Appendix \ref{ap1}.

Taking into account that $t\in T$ and $s\in S$ to solve system
\eqref{r1}, we see that this can be reduced to solve the system $H_1=0$, $H_2=0$,
because $s (2s-1) (s^2+2s-1)\ne 0.$ Hence, to find the real
solutions of the systems $h_1=0,$ $  h_2=0,$ is equivalent to find
the real solutions, $(s,t)\in S\times T,$ of the polynomial system
$H_1=0,$ $H_2=0.$

To study the above planar systems of equations we will use a mixture
of the classical approach applying resultants together with simple
inequalities in the original variables $x_3$ and $y_5.$ For our
problem this approach is very suitable because of the simplicity of
the equation $h_2=0,$ given in~\eqref{eq:h2}.

We start with the polynomials system $H_1(s,t)=0,$ $H_2(s,t)=0.$
Recall that if $(\hat s,\hat t)$ is one of its solutions (real or
complex), then $t=\hat t$ must be a zero of the one variable polynomial
\[
P(t)=\operatorname{Res}_s(H_1,H_2),
\]
where $\operatorname{Res}_s(\cdot,\cdot)$ denotes the resultant of
two polynomials with respect the variable $s,$ see for instance
\cite{Sturm}. After some computations (implemented for instance in
Maple or Mathematica) we get that
\begin{equation}\label{eq:pt}
P(t)=(1 + t^2)^6 p_4(t)q_4(t)p_{120}(t)p_{132}(t),
\end{equation}
where $p_4(t)=1 + 4 t - 14 t^2 + 4 t^3 + t^4,$ $q_4(t)=1 - 4 t - 14
t^2 - 4 t^3 + t^4$ and $p_k$ denotes a polynomial with integer
coefficients of degree $k$ that we do not detail. We only remark that precisely $p_{120}$ is the polynomial that gives rise to the polynomial $R_{60},$ detailed in Appendix \ref{ap2}  and that gives rise to the values $t^*$ and $u^*$ that appear in Remark \ref{re:1}. Hence, by
computing the Sturm sequences of each of the four polynomials,
$p_4,q_4,p_{120}$ and $p_{132},$ we get that they have respectively,
4,4, 28 and 32 real roots and, moreover, that all them are simple.
Furthermore, since we are only interested on the roots $t\in
T\sim(0.127,1),$ we consider a slightly  bigger interval $T\subset
T'=(3/25,100),$ with rational endpoints. Again, the corresponding
Sturm sequences allow to prove that their number of roots in $T'$
are $1,1,7$ and $9,$ respectively. We will denote them by $t_1;$
$t_2;$ $t_3,\ldots,t_9$ and $t_{10},\ldots t_{18},$ where for each
$p_k$ the roots are ordered. Their approximated ordered value are
$$
\begin{array}{lllll}
t_{10}\approx 0.1278827,& t_{3}\approx 0.1296657,& t_{11}\approx
0.1318307,& t_{12}\approx 0.1535285,& t_{2}\approx 0.1583844,\\
t_{13}\approx 0.1690804,& t_{4}\approx 0.1818971,& t_{5}\approx
0.1871837,& t_{14}\approx 0.4693713,& t_{1}\approx 0.5095254,\\
t_{15}\approx 0.5490528,& t_{16}\approx 0.5930556,& t_{6}\approx
0.7095411,& t_{7}\approx 0.7332148,& t_{8}\approx 0.9432977,\\
t_{17}\approx 0.9681690,& t_{9}\approx 0.9958185,& t_{18}\approx
0.9962499.&&
\end{array}
$$
In fact, the roots $t_1$ and $t_2$ are  $t_1=-1+\sqrt5-\sqrt{5-2\sqrt5} \approx 0.5095254$ and $t_2=1+\sqrt5-\sqrt{5+2\sqrt5}\approx0.1583844,$
while each of the other sixteen roots can be obtained, again using
the Sturm sequences, with any desired error.

Therefore, each $t_j, j=3,\ldots,18,$  can be bounded by
$\underline{t}_j<t_j<\overline{t}_j,$ with
$\underline{t}_j,\overline{t}_j\in\mathbb{Q}$
and $0<\overline{t}_j-\underline{t}_j$ as small, as desired. For
each of these values of $t_j$ we can use that the function $t\to
(1-t^2)/(4t)$ is decreasing in $(0,1)$ and that $y_5=(1-t^2)/(4t)$
to get that if $(1-t_j^2)/(4t_j)=y_5(j),$ then
\[
\underline y_5(j)=\frac{1-\overline t_j^2}{4\overline t_j}<y_5(j)<
\frac{1-\underline t_j^2}{4\underline t_j}=\overline y_5(j),
\quad\mbox{with}\quad \underline y_5(j),\overline
y_5(j)\in\mathbb{Q}.
\]

Recall that in \eqref{eq:psi}  we have introduced the functions
$$
\Psi^\pm(y)=\frac14 \pm\frac{y}2\sqrt{\frac{15-4y^2}{1+4y^2}}.
$$
By simple derivation we get that $\Psi^+$  (resp.   $\Psi^-$ )
is increasing (resp. decreasing) for $y\in(0,\sqrt3/2)$ and decreasing (resp. increasing) for  $y\in(\sqrt3/2,1).$ Hence, if we define
$
x_3^\pm=\Psi^\pm( y_5)
$
it holds that
\begin{itemize}
    \item  If $j$ is such that $y_5(j)<\sqrt3/2$ then
    \[
    \Psi^+(\underline{y}_5(j))< x_3^+(j) <\Psi^+(\overline{y}_5(j)),\quad     \Psi^-(\overline{y}_5(j))< x_3^-(j)  <\Psi^-(\underline{y}_5(j)).
    \]
    \item  If $j$ is such that $y_5(j)>\sqrt3/2$ then
    \[
    \Psi^+(\overline{y}_5(j))< x_3^+(j)  <\Psi^+(\underline{y}_5(j)),\quad    \Psi^-(\underline{y}_5(j))< x_3^-(j) <\Psi^-(\overline{y}_5(j)).
    \]
    \end{itemize}
From these inequalities it is easy to find rational values $\underline x^\pm_3(j)$ and $\overline x^\pm_3(j)$ such that
\[
 \underline x^\pm_3(j)<x^\pm_3(j)< \overline x^\pm_3(j),
\]
and with $0< \overline x^\pm_3(j)- \underline x^\pm_3(j)$ as small
as desired. Finally, the functions that define $E_1,E_2$ and $E_3,$
given by $\psi_1(x)=\sqrt{1+2x},$ $\psi_2(x)=\sqrt{3+4x-4x^2}$ and
$\Psi_3(y)=(1+4y^2)^{3/2}$ are increasing,  increasing for
$x\in(0,\sqrt3/3)$ and decreasing for $x\in(\sqrt3/3,1),$ and
increasing, respectively. Similarly that for $x^\pm(j)$ we can find
rational bounds, $\underline E_j$ and $\overline E_j,$ $j=1,2,3,$ as sharp as desired and satisfying
\[
\underline E_1<E_1<\overline E_1, \quad \underline E_2<E_2<\overline E_2, \quad \underline E_3<E_3<\overline E_3.
\]

In short, we have found sharp rational upper and lower bound of any
possible solution $y_5=y_5(j),$ $x_3^\pm(j),$ corresponding to each
$t=t_j,j=1,2,\ldots,18.$ These rational bounds give also rational
bounds for $E_1,E_2$ and $E_3.$  Gluing these bounds we prove that
some candidates to be solutions of our system can be discarded. That
other candidates are actual solution of our problem can be easily
proved from Bolzano's theorem.
 We detail two examples, one
of each type and skip the computations for all the rest.

A suitable way is to write $h_1$ in function of the  variables $x_3,y_5,E_1,E_2$ and $E_3$, namely
\begin{align*}
      h_1=&-128{{  x_3}}^{3}{{  y_5}}^{2} \left( 4{{  x_3}}^{2}+4{  x_3
    }+1 \right) {  E_ 1}+4{{  x_3}}^{2}{{  y_5}}^{2} \left( 4{{
            x_3}}^{2}+4{  x_3}+1 \right) {  E_ 3}\\&+32{{  y_5}}^{2} \left( 32
    {{  x_3}}^{6}+16{{  x_3}}^{5}-8{{  x_3}}^{4}-4{{  x_3}}^{3}+
    4{{  x_3}}^{2}+4{  x_3}+1 \right) {{  E_ 1}}^{2}+64{{  x_3}}^
    {3}{  y_5} \left( 1+2{  x_3} \right) {  E_ 1}{  E_ 2}\\&-8{{
              x_3}}^{2}{{  y_5}}^{2} \left( 4{{  x_3}}^{2}+4{  x_3}+1
    \right) {  E_ 1}{  E_ 3}-4{{  x_3}}^{2}{  y_5} \left( 1+2{
          x_3} \right) {  E_ 3}{  E_ 2}\\&+16{  y_5} \left( 16{{  x_3}
    }^{5}+16{{  x_3}}^{4}+4{{  x_3}}^{3}-4{{  x_3}}^{2}-4{
        x_3}-1 \right) {{  E_ 1}}^{2}{  E_ 2}\\&-4{{  y_5}}^{2} \left( 16{{
              x_3}}^{6}-16{{  x_3}}^{5}-32{{  x_3}}^{4}-12{{  x_3}}^{3}+
    3{{  x_3}}^{2}+4{  x_3}+1 \right) {{  E_ 1}}^{2}{  E_ 3}\\&+8{{
              x_3}}^{2}{  y_5} \left( 2{{  x_3}}^{2}+3{  x_3}+1 \right)
    {  E_ 1}{  E_ 2}{  E_ 3}+{{  E_ 2}}^{2}{  E_ 3}{{  x_3}}^{2}\\&-2
    {  y_5} \left( 32{{  x_3}}^{5}+40{{  x_3}}^{4}+16{{  x_3
    }}^{3}-2{{  x_3}}^{2}-4{  x_3}-1 \right) {{  E_ 1}}^{2}{  E_ 2}
    {  E_ 3}\\&-2{{  x_3}}^{2} \left( 1+2{  x_3} \right) {  E_ 1}{{
              E_ 2}}^{2}{  E_ 3}+{{ x_3}}^{2} \left( 4{{x_3}}^{2}+4{
          x_3}+1 \right) {{E_1}}^{2}{{E_2}}^{2}{E_3}.
    \end{align*}
We remark that the addends of $h_1$ can be bounded using the inequalities given above.
For instance
\[
{{  \underline E_ 2}}^{2}{ \underline E_ 3}{{ \underline
x_3}}^{2}<{{  E_ 2}}^{2}{  E_ 3}{{ x_3}}^{2}<{{\overline  E_
2}}^{2}{ \overline E_ 3}{{ \overline x_3}}^{2},
\]
and
\[
-4{{\overline  x_3}}^{2}{ \overline y_5} \left( 1+2{
        \overline  x_3} \right) {\overline  E_ 3}{\overline  E_ 2}<-4{{  x_3}}^{2}{  y_5} \left( 1+2{
          x_3} \right) {  E_ 3}{  E_ 2}<-4{{\underline  x_3}}^{2}{ \underline y_5} \left( 1+2{
      \underline    x_3} \right) {\underline  E_ 3}{\underline  E_ 2}.
\]

Let us prove for instance that the value $y_5=y_5(5)$ corresponding
to $t=t_5\approx 0.1871837,$ together with $x_3=x_3^+(5)=\Psi^+(
y_5(5))$ does not provide a solution of our system. From the Sturm
sequence we get that
\[
\frac{1871}{10000}<t_5<\frac{1872}{10000}.
\]
From these inequalities and all the above considerations we obtain
that
\begin{align*}
1.2886701\approx &\frac{94234}{73125}<
y_5<\frac{96499359}{74840000}\approx 1.2894088,\quad
\frac{9235}{10000}<x_3^+<\frac{9238}{10000},\\
&\frac{16873}{10000}<E_1< \frac{16873}{10000},\quad
\frac{18118}{10000}<E_2<\frac{18115}{10000},\quad
\frac{21128}{1000}<E_3<\frac{21161}{1000}.
\end{align*}
By using all these inequalities we obtain that the corresponding
value of $h_1>242$ and the system has no solution  for the value of
$y_5,$ and its corresponding $x_3^+(5)$ associated to $t=t_5.$

On the other hand, let us prove that the value $y_5=y_5(7)$
corresponding to  $t=t_7=t^*\approx 0.7332148$ together with
$x_3=x_3^+(7)=\Psi^+( y_5(7))$ does provide an actual solution.
This is a simple consequence of Bolzano's theorem, because if we
denote as $h_1(\tau)$ the value of the expression of $h_1$ when all
the values $y_5,x_5^+,E_1,E_2$ and $E_3$ are obtained when $t=\tau$
we get for instance that
\[
h_1\Big(\frac{7332}{10000}\Big)h_1\Big(\frac{7333}{10000}\Big)<0.
\]

We carry out similar computations for $y_5=y_5(j),$ corresponding to
$t=t_j,$ and $x_3=x_3^\pm(j)=\Psi^\pm (y_5(j))$.
We conclude that among the 36 candidates that could be a solution of the system
$ h_1(x_3,y_5)=0$, $h_2(x_3,y_5)=0$ the only couples $(x_3,y_5)$ that do solve it are:
\begin{enumerate}[(I)]
\item $(x_3^+(2),y_5(2))=\left(\dfrac{1+\sqrt5}4,  \dfrac{\sqrt{5+2\sqrt5}}2\right)
\approx(0.8090170,1.5388418)$  corresponding to $t=t_2,$
\item $(x_3^+(7),y_5(7))\approx(0.5402091,0.1576605)$ corresponding to $t=t_7=t^*,$
\item $(x_3^+(9),y_5(9))\approx(0.2540572,0.0020951)$ corresponding to $t=t_9,$
\item $(x_3^-(4),y_5(4))\approx(-0.4091526,1.3289291)$ corresponding to $t=t_4,$
\item $(x_3^-(14),y_5(14))\approx(-0.3542470,0.4152845)$ corresponding to $t=t_{14},$
\item $(x_3^-(18),y_5(18))\approx(0.2463622,0.0018786)$ corresponding to $t=t_{18}$.
\end{enumerate}
Clearly, solutions  in items (IV) and (V) can be discarded
because the corresponding values of $x_3$ are negative. The solutions given in items (III) and (VI) are not good either, since both options result in negative $m_5$ values.
In short the central configurations are:
\begin{enumerate}
\item [(I)] $(x_3,y_3)=\left(\dfrac{1+\sqrt5}4,
\dfrac{\sqrt{10+2\sqrt5}}4\right)$ and
$y_5=\dfrac{\sqrt{5+2\sqrt5}}2,$ with masses $m_j=\frac1 5,$ for all
$j.$
\item [(II)] $(x_3,y_3)=(x_3^+(7),y_3^+(7))\approx(0.5402091,0.9991913)$ and $y_5=y_5(7)\approx0.1576605$  with
masses $m_1=m_2\approx 0.0922539,$ $m_3=m_4\approx0.3860949$ and
$m_5 \approx 0.04330243.$
\end{enumerate}
The values of $y_3$ and  $m_j$ are obtained from Proposition
\ref{p1} and  the expressions \eqref{eq:m1m3} and \eqref{eq:lambda}.

The first solution provides the regular pentagon of Figure 2(a) as a
convex central configuration of the $5$-body problem with masses
equal to $1/5$. While the second one provides the equilateral
concave pentagon of Figure 2(b) as a concave central configuration
of the $5$-body problem with the masses given in the statement (b)
of the theorem. This completes the proof of Theorem \ref{t1}.

\subsection{Alternative approaches to solve system $  h_1=0, h_2=0$} In this section
 we comment about to alternative approaches two solve this system and its equivalent one $H_1=0, H_2=0.$

A first one consists on computing the Gr\"obner basis of the two
polynomials $H_1$ and $H_2$ with respect to the two variables $s$
and $t.$ Doing this  we obtain three polynomial equations whose
common solutions are also the solutions of system $H_1=0$, $H_2=0$.
We do not provide explicitly these three polynomials, but only
comment that they are huge and their expressions need many pages.
The first one essentially coincides with $P(t)$ given in
\eqref{eq:pt}. The second one $ P_2(s,t)= (1 + t^2)^2 p_{259}(s,t).
$ depends on both variables $s$ and $t$, and $P_2$ is linear in the
variable $s$. Consequently each root $t=t_j$  of the polynomial $P$
provides a single value of $s$ from $P_2(s,t_j)=0,$ say $s=s_j.$
Then we only need to keep  the $j'$s such that $s_j\in S.$ Finally,
the third polynomial $P_3(s,t)$ of the Gr\"obner basis has degree
$262$ but it is only cubic in the variable $s$. By keeping only the
values $s_j\in S$ that also satisfy $P_3(s_j,t_j)=0$ we arrive to
the actual solutions. We have used  our approach instead of this one
because it is not easy to check analytically all the above facts
because only  two of the eighteen roots of $P$ are known
analytically. Moreover, we prefer our point of view because the
computation of a resultant is simple and self contained while the
computation of a Gr\"obner basis is implemented in the computer
algebra systems but the user has no control on  what the algorithm
is doing.

A second alternative approach would consist on computing also
$Q(s)=\operatorname{Res}_t(H_1,H_2).$ In this case we arrive to
\begin{equation*}
    Q(s)=q_2^6(s) q_4^2(s)q_{120}(s)q_{132}(s),
\end{equation*}
for some polynomials $q_k$ of degree $k,$ where here $q_4$ is
different to  the one given in \eqref{eq:pt}. Their respective
number of real roots are $0,4,28$ and $32.$ Moreover, only $1,0,4$
and $6$ of them are in $S.$ Call them $s_m,$ $m=1,2,\ldots,11.$
Hence all the possible solutions of system $H_1=0,H_2=0$ in $S\times
T$ are given by $11\times 18$ values $(s_m,t_j)$. Then, a discarding
process, similar to the one done in our proof of Theorem~\ref{t1}
can be done. On the other hand, a proof that the non discarded
candidates to be solutions are actual solutions can be done for
instance by using the nice Poincar\'{e}-Miranda theorem. See for
instance~\cite{GasLloMan} to have more details of how to utilize this
approach. We have not used it in our work because the expression of
$h_2$ is much simpler that the one of $Q.$

\appendix

\section{The expression of $H_1$}\label{ap1}
The polynomial $H_1$ in \eqref{r1} writes as
\begin{align*}
    H_1(s,t)=
    R_{0}(s)\big(t^{10}+1)-R_1(s)\big(t^9+4t^7-t\big)+R_2(s)\big(t^8+t^2\big)+R_3(s)\big(t^7+t^3\big)+\sum_{j=4}^6
    R_j(s)t^j,
\end{align*}
where
\begin{align*}
R_0=&{s}^{2} \left( 2 s-1 \right) ^{ 2}\Big(
140137001{s}^{20}-473336800{s}^{19}-77771662{s}^{
18}+3288160224{s}^{17}\\&-8637659443{s}^{16}+12537556864{s}^{15}-
12225124968{s}^{14}+8691543680{s}^{13}\\&-4785798270{s}^{12}+
2180211392{s}^{11}-880204628{s}^{10}+324883264{s}^{9}-105646862
{s}^{8}\\&+28078976{s}^{7}-5656040{s}^{6}+796288{s}^{5}-67435{s
}^{4}+1568{s}^{3}+306{s}^{2}-32s+1 \Big),\\
R_1=&2s \left( s-1 \right)  \left( 3s-1 \right)  \left( 2s-1 \right)
\Big( 402088273{s}^{20}-1176961940{s}^{19}+169440330{s}^{18}\\&+
4244422908{s}^{17}-9214008723{s}^{16}+10491664368{s}^{15}-
7879610248{s}^{14}+4350977648{s}^{13}\\&-1966298574{s}^{12}+
810673640{s}^{11}-313304452{s}^{10}+104139528{s}^{9}-26338798{
s}^{8}\\&+4329392{s}^{7}-227208{s}^{6}-97872{s}^{5}+33277{s}^{4}-
5620{s}^{3}+586{s}^{2}-36s+1 \Big),\\
R_2=&2610735845{s}^{24}-15402964740{s}^{23}+36242909513{s}^{22}-
34656361080{s}^{21}-19279444736{s}^{20}\\&+100341955724{s}^{19}-
143714126267{s}^{18}+121841849904{s}^{17}-68503362257{s}^{16}\\&+
27259846072{s}^{15}-9818466526{s}^{14}+5237546480{s}^{13}-
3390994016{s}^{12}+1681689656{s}^{11}\\&-535976190{s}^{10}+78574400
{s}^{9}+17032379{s}^{8}-14081204{s}^{7}+4575253{s}^{6}\\&-963768
{s}^{5}+142992{s}^{4}-15044{s}^{3}+1081{s}^{2}-48s+1,\\
R_3=&4s \left( 2s-1 \right)  \Big( 4011437555{s}^{22}-17997365760{
s}^{21}+47815887103{s}^{20}-119864245704{s}^{19}\\&+262594976801{s}
^{18}-427430598888{s}^{17}+498867338773{s}^{16}-424920266400{s}^
{15}\\&+272105947982{s}^{14}-136189975712{s}^{13}+56438427910{s}^{
12}-20882782960{s}^{11}\\&+7255673746{s}^{10}-2308064816{s}^{9}+
623396538{s}^{8}-133508384{s}^{7}+21674927{s}^{6}\\&-2595360{s}^{
5}+227627{s}^{4}-15240{s}^{3}+861{s}^{2}-40s+1 \Big),\\
R_4=&-128824963717{s}^{24}+476911508392{s}^{23}+165019719846{s}^{22}-
4572358588576{s}^{21}\\&+13580634656322{s}^{20}-22953475876792{s}^{
19}+26679405170926{s}^{18}-23043609894160{s}^{17}\\&+15476921278317
{s}^{16}-8354767360048{s}^{15}+3726249539948{s}^{14}-1406580864896
{s}^{13}\\&+458981618588{s}^{12}-132125593264{s}^{11}+34342034540
{s}^{10}-8246853600{s}^{9}+1837269525{s}^{8}\\&-367787384{s}^{7}+
62442990{s}^{6}-8479328{s}^{5}+879394{s}^{4}-67480{s}^{3}+3750
{s}^{2}-144s+3,\\
\end{align*}

\begin{align*}
R_5=&-1024{s}^{3} \left( 2s-1 \right) ^{3} \Big(
5478853{s}^{18}-19634676{s}^{17}+61412857{
s}^{16}-189035488{s}^{15}\\&+395309060{s}^{14}-534006384{s}^{13}+
487648452{s}^{12}-316317344{s}^{11}+154130630{s}^{10}\\&-61395832
{s}^{9}+22312398{s}^{8}-7766560{s}^{7}+2406420{s}^{6}-594416{s
}^{5}\\&+108148{s}^{4}-13664{s}^{3}+1117{s}^{2}-52s+1
\Big),\\
R_6=&138883898747{s}^{24}-544790344792{s}^{23}+35213761702{s}^{22}+
4234516514656{s}^{21}\\&-13218500821438{s}^{20}+22675637356616{s}^{
19}-26464806060818{s}^{18}+22812618679024{s}^{17}\\&-15226596524179
{s}^{16}+8152770613200{s}^{15}-3614885715604{s}^{14}+1368000514176
{s}^{13}\\&-453304331364{s}^{12}+133617865040{s}^{11}-35152022420
{s}^{10}+8179417120{s}^{9}-1643965931{s}^{8}\\&+272577160{s}^{7}-
34266642{s}^{6}+2694560{s}^{5}-25822{s}^{4}-22424{s}^{3}+2726
{s}^{2}-144s+3
\end{align*}

\vspace{-5.5cm}

\section{The polynomial $R_{60}$}\label{ap2}
The polynomial $p_{120}(t)$ is reciprocal, that is
$p_{120}(t)-t^{120}p_{120}( 1/t)\equiv0.$
Notice that if $\hat t\,$ is one of its roots, $1/\hat t$ is another one. Hence there is a standard trick to ``reduce" its degree to the half. Consider the numerator of $t+1/t=u,$ that is $t^2+tu-1$ and compute the resultant  between it and $p_{120}(t).$ We obtain that
\[\operatorname{Res}_t\big(p_{120}(t),t^2-ut+1\big)= (R_{60}(u))^2,\]
where $R_{60}$ is the   polynomial of degree $120/2=60,$
\begin{align*}
&{u}^{60}+4{u}^{59}-480{u}^{58}-2368{u}^{57}+102656{u}^{56}+
661504{u}^{55}-12378112{u}^{54}-114393088{u}^{53}\\&+813367296{u}
^{52}+13487570944{u}^{51}-6779043840{u}^{50}-1119258411008{u}^{
    49}\\&-4587041849344{u}^{48}+63761809408000{u}^{47}+556458915659776
{u}^{46}-2085902406909952{u}^{45}\\&-38793866899357696{u}^{44}-
14062083704356864{u}^{43}+1871701598900584448{u}^{42}\\&+
6613561341561012224{u}^{41}-63439436081954553856{u}^{40}-
468879147034843545600{u}^{39}\\&+1328797291115756650496{u}^{38}+
20993318838230010822656{u}^{37}\\&-2305632465931114381312{u}^{36}-
680737782622703312699392{u}^{35}\\&-1102380299141445066948608{u}^{34}
+16433802069777919820955648{u}^{33}\\&+51705311821812917239545856{u}^
{32}-287454839290286351637807104{u}^{31}\\&-
1467660179422186654016733184{u}^{30}+3159833116868066015124127744{
    u}^{29}\\&+30235824601650376596023934976{u}^{28}-
3647821652057127970278473728{u}^{27}\\&-478786881979744683498227630080
{u}^{26}-692146369158869263721793847296{u}^{25}\\&+
6031919651486804311277903544320{u}^{24}+
17009050671171345412955758919680{u}^{23}\\&-
62657388272839173181632407404544{u}^{22}-
249363510060237095878177991950336{u}^{21}\\&+
558394430418773435280130962358272{u}^{20}+
2603575445697403174765988761567232{u}^{19}\\&-
4366163704989486187475739888058368{u}^{18}-
20183899906055516645882016026329088{u}^{17}\\&+
29417136038642440588505535383339008{u}^{16}+
117121419223047038546206624994820096{u}^{15}\\&-
162680894447793178205603154016337920{u}^{14}-
505682422731087446635760756082081792{u}^{13}\\&+
700532798421762302278975021920550912{u}^{12}+
1597276835067625721595839817342517248{u}^{11}\\&-
2240245382048536583959836766075092992{u}^{10}-
3574566797945169482598857579834638336{u}^{9}\\&+
5053220375882588124433490027151884288{u}^{8}+
5336383096359169095222217602355953664{u}^{7}\\&-
7423686433490169801891477126702432256{u}^{6}-
4673067172681344865677446696298086400{u}^{5}\\&+
5934795309305307979410357304298569728{u}^{4}+
1599227432428726909587392869399789568{u}^{3}\\&-
996920996838686904677855295210258432{u}^{2}+
332306998946228968225951765070086144u\\&-
1329227995784915872903807060280344576.
\end{align*}
Then all the roots of $p_{120}$ can be obtained simply by computing the roots of $R_{60}$ and then for each one of them, say $u=\hat u,$ two roots of $p_{120}$ are given by the solutions of the quadratic equation $t^2-\hat u t+1=0.$

By computing its Sturm sequence we get that $R_{60}$ has exactly $14$ real roots, all them simple (of course half the number the real roots of $p_{120}$). One of them is $u=u^*\approx  2.0970716051$ and the value $t=t^*$ appearing in Remark \ref{re:1} is the smallest root of $t^2-u^*t+1=0.$

\section*{Acknowledgements}
The first author  is partially supported by  the grant
Sistemas Hamiltonianos, Mec\'anica y Geometr\'{\i}a from the PAPDI2021 CBI-UAMI.
The last two authors are partially by the Agencia Estatal de Investigaci\'on grant PID2019-104658GB-I00,  and the H2020 European Research Council grant MSCA-RISE-2017-777911.

\doublespacing
\bibliographystyle{abbrvnat}

\bibliography{referen5body}

\end{document}